\documentclass[12pt,a4wide]{amsart}
\usepackage{amsmath,amsthm,amsfonts,amssymb,graphicx}

\usepackage[all]{xy}

\newtheorem{theorem}{Theorem}[section]
\newtheorem{lemma}[theorem]{Lemma}
\newtheorem{definition}[theorem]{Definition}

\newtheorem{proposition}[theorem]{Proposition}

\newcommand\Hom{\operatorname{\mathcal{H}om}}
\newcommand\ho{\operatorname{Hom}}

\newcommand\im{\operatorname{Im}}

\newcommand\ex{\operatorname{Ext}}

\newcommand\Pic{\operatorname{Pic}}

\newcommand\rk{\operatorname{r}}

\newcommand\dg{\operatorname{d}}

\title[Relative Jacobians of elliptic fibrations]
{Relative Jacobians of elliptic fibrations with reducible fibers}
\author{Ana Cristina L\'opez Mart\'{\i}n }
\email{anacris@usal.es}
\address{Departamento de Matem\'aticas, Universidad de Salamanca,
Plaza de la Merced 1-4, 37008 Salamanca, Spain}
\date{\today}
\thanks {e-mail: anacris@usal.es. \\ This research was partly supported by the Spanish-Italian cooperation project
HI00-141, by the research projects BFM2000-1315 and SA009/01 and
by EAGER Contract no. HPRN-CT-2000-00099} 

\begin{document}
\maketitle
\begin{abstract} We prove that if $X$ and  $S$ are smooth varieties and $f\colon X\to S$ is an
elliptic fibration with singular fibers curves of types I$_N$ with
$N\geq 1$, II, III and IV, then the relative Jacobian
$\hat{f}\colon \overline{M}_{X/S}\to S$ of $f$, defined as the
relative moduli space of semistable pure dimension one sheaves of
rank 1 and degree 0 on the fibers of $f$, is an elliptic fibration
such that all its fibers are irreducible. This extends known
results when fibers are integral or of type I$_2$.
\end{abstract}
\section{Introduction}

Let $f\colon X\to S$ be an elliptic fibration, that is, a proper
flat morphism of schemes whose fibers are Gorenstein curves of
arithmetic genus 1. If we have a relative ample sheaf
$\mathcal{O}_X(1)$ on $X$, the relative Jacobian of $f$ is defined
as the Simpson moduli space $\hat{f}\colon \overline{M}_{X/S}\to
S$ of semistable pure dimension one sheaves of rank 1 and degree 0
on the fibers of $f$ with respect to the polarization induced by
$\mathcal{O}_X(1)$. Since every torsion free rank 1 sheaf on an
integral curve is stable, when $f\colon X \to S$ is an integral
fibration, that is, all its fibers are geometrically integral
curves, its relative Jacobian is simply the relative moduli space
of torsion free rank 1 sheaves of relative degree 0. In this case,
it is known that for any closed point $s\in S$ we have an
isomorphism $(\overline{M}_{X/S})_s\simeq X_s$ between the fibers
of both fibrations over $s$. Then this relative Jacobian is again
an integral elliptic fibration. Furthermore if the fibration
$f\colon X\to S$ has a section, $\overline{M}_{X/S}$ is not only
locally but globally isomorphic to $X$ over $S$ (see \cite{HM}).

On the other hand the geometry of such relative Jacobian
$\hat{f}\colon \overline{M}_{X/S}\to S$ is not known for a general
elliptic fibration. In \cite{Cal2} we find some examples showing
that when $f\colon X\to S$ has reducible fibers (exactly $I_2$
fibers), it is no longer true that the two fibrations $f$ and
$\hat{f}$ have isomorphic fibers. In fact, C\u ald\u araru proves
that if $X_s$ is a fiber of type $I_2$, then the fiber
$(\overline{M}_{X/S})_s$ is isomorphic to a rational curve  with
one node.

One should point out that elliptic fibrations have been used in
string theory, notably in connection with mirror symmetry on
Calabi-Yau manifolds and D-branes. Some of the classic examples of
families of Calabi-Yau manifolds for which there is a
Greene-Plesser \cite{GP90} description of the mirror family
produced by Candelas and others, are elliptic fibrations
\cite{Can91}. Moreover there is a relative Fourier-Mukai transform
for most elliptic fibrations (\cite{BBHM, HM}) that can be
understood in terms of duality in string theory (\cite{DoTani,
ACHY} or D-brane theory. The latter application is due to the
interpretation of B-type D-brane states as objects of the derived
category $D(X)$ of coherent sheaves \cite{Kon95, As, Dou} and  to
Kontsevich's  homological mirror symmetry proposal \cite{Kon95}
that gives an equivalence between  $D(X)$ and the Fukaya category
\cite{Fuk93}. According to that proposal, the monodromies around
special points of the known-models of local moduli spaces of
Lagrangian submanifolds should correspond to Fourier-Mukai
transforms \cite{ACHY, Ho}; this explains the importance of
elliptic Calabi-Yau manifolds in string theory.

The aim of this paper is to study the structure of the relative
Jacobian of an elliptic fibration $f\colon X\to S$ such that $X$
and $S$ are smooth projective varieties and the singular curves
appearing as fibers of $f$ are:

(I$_1$) : A rational curve with one node.

(II) : A rational curve with one cusp.

(III) : $C_1\cup C_2$ where $C_1$ and $C_2$ are rational smooth
curves with $C_1\cdot C_2=2p$.

(IV) : $C_1\cup C_2\cup C_3$, where $C_1, C_2, C_3$ are rational
smooth curves and $C_1\cdot C_2=C_2\cdot C_3=C_3\cdot C_1=p$

(I$_N$) : $C_1\cup C_2\cup \hdots \cup C_N$,  where $C_i$,
$i=1,\hdots, N$, are rational smooth curves and $C_1\cdot
C_2=C_2\cdot C_3=\hdots=C_{N-1}\cdot C_N=C_N\cdot C_1=1$ if $N>2$
and $C_1\cdot C_2=p_1+p_2$ if $N=2$.

If $f\colon X\to S$ is an elliptic fibration of this type, we
prove that the moduli space of semistable pure dimension one
sheaves of rank 1 and degree 0 on a fiber $X_s$ is isomorphic to a
smooth elliptic curve when $X_s$ is smooth, to a rational curve
with one node when $X_s$ is I$_N$, $N\geq 1$ and to a rational
curve with one cusp when $X_s$ is II, III or IV. The result is
then that the relative Jacobian $\hat{f}\colon
\overline{M}_{X/S}\to S$ of $f$ is an integral elliptic fibration
(Theorem \ref{teorema}). In particular we deduce that, although
$\overline{M}_{X/S}$ is irreducible, if the fibration $f$ has
reducible fibers, it cannot be isomorphic to $X$, even assuming
that $f$ has a section. For instance, if the variety $X$ has
dimension 2 or 3, then the results about contractibility of curves
on smooth surfaces \cite{A} and on smooth threefolds \cite{Z}
allow us to  to ensure the existence of singular points in the
moduli space $\overline{M}_{X/S}$ that correspond to strictly
semistable sheaves on the fibers of $f$. Then in these cases
$\overline{M}_{X/S}$ is not isomorphic to the original variety
$X$.

Since by Kodaira's work \cite{K}, we known that every elliptic
surface $f\colon X\to S$ with reduced fibers is a fibration of
this type, the theorem gives us in particular the structure of the
relative Jacobian of any reduced elliptic surface.

I wish to express my gratitude to Prof. Hern\'andez Ruip\'erez for
his invaluable help and for his constant encouragement. I also
thank  A. C\u ald\u araru for many useful suggestions.

\section{Preliminares} All the schemes considered in this paper
are of finite type over an algebraically closed field $\kappa$ of
characteristic zero and all the sheaves are coherent.

\subsection{Some elliptic fibrations with reducible fibers}
 Let $f\colon X\to S$ be
an elliptic fibration. By this we mean a proper flat morphism  of
schemes whose fibers are geometrically connected Gorenstein curves
of arithmetic genus 1. We denote by $X_s$ the fibre of $f$ over
$s\in S$ and by $\Sigma(f)$ the {\it discriminant locus} of $f$,
that is, the closed subset of points $s\in S$ such that $X_s$ is
not a smooth curve.

There are two important cases where the curves that can occur as
singular fibers of an elliptic fibration $f\colon X\to S$ are
known. When $X$ is an elliptically fibred surface over a smooth
curve $S$, the singular fibers of $f$ are classified by Kodaira
\cite{K}. And if $f\colon X\to S$ is an elliptic threefold with
$X$ and $S$ smooth projective varieties and the map $f$ has a
section, Miranda \cite{M} studies the kinds of degenerated fibers
of $f$ that can appear. These two works allow us to ensure the
existence of elliptic fibrations as in the following

\begin{definition}\label{d:fibracion} An elliptic fibration of type ($*$) is an
elliptic fibration $f\colon X\to S$, with $X$ and $S$ smooth
projective varieties, without multiple fibers and such that if
$s\in \Sigma(f)$, the fiber $X_s$ is one of the following curves:

(I$_1$) : A rational curve with one node.

(II) : A rational curve with one cusp.

(III) : $X_s=C_1\cup C_2$ where $C_1$ and $C_2$ are rational
smooth curves with $C_1\cdot C_2=2p$.

(IV) : $X_s=C_1\cup C_2\cup C_3$, where $C_1, C_2, C_3$ are
rational smooth curves and $C_1\cdot C_2=C_2\cdot C_3=C_3\cdot
C_1=p$

(I$_N$) : $X_s=C_1\cup C_2\cup \hdots \cup C_N$,  where $C_i$,
$i=1,\hdots, N$, are rational smooth curves and $C_1\cdot
C_2=C_2\cdot C_3=\hdots=C_{N-1}\cdot C_N=C_N\cdot C_1=1$ if $N>2$
and $C_1\cdot C_2=p_1+p_2$ if $N=2$.
\end{definition}

If $S$ is a smooth curve and $f\colon X\to S$ is an elliptic
surface such that all fibers are reduced, by Kodaira's
classification, every singular fiber $X_s$ of $f$ is one of the
above list. Then any smooth elliptic surface over $S$ with reduced
fibers is a elliptic fibration of type ($*$).

If $f_0\colon X_0\to S_0$ is an elliptic fibration with $X_0$ and
$S_0$ varieties of dimensions 3 and 2 respectively and $f_0$ has a
section, Miranda constructs in \cite{M} a flat model $f\colon X\to
S$ with $X$ and $S$ smooth and such that the discriminant locus
$\Sigma(f)$ is a curve with at worst ordinary double points as
singularities. He proves that at a smooth point $s\in \Sigma(f)$
the singular fiber $X_s$ is one on Kodaira's list and that the
type of fiber is constant on the irreducible components of
$\Sigma(f)^{\text{smooth}}$. Moreover he determines, case by case,
the type of singular fiber over the collision points of
$\Sigma(f)$. The results of Miranda imply that if the fibers of
$f$ over the smooth points of $\Sigma(f)$ are reduced and all
collisions are of type $I_{N_1}+I_{N_2}$, then $f\colon X\to S$ is
an elliptic fibration of type ($*$).

\subsection{The Jacobian of a projective curve}
Let $C$ be a projective curve. Let ${\mathcal L}$ be an ample
invertible sheaf on $C$, let $H$ be the associated polarization
and let $h$ denote the degree of $H$.

A sheaf $F$ on $C$ is {\it pure dimension one} if the support of
any nonzero subsheaf of $F$ is of dimension one. The (polarized)
{\it rank} and {\it degree} with respect to $H$ of $F$ are the
rational numbers $\rk_H(F)$ and $\dg_H(F)$ determined by the
Hilbert polynomial
$$P(F,n,H)=\chi(F\otimes\mathcal{O}_C(nH))=h\rk_H(F)n+\dg_H(F)+\rk_H(F)\chi(\mathcal{O}_C).$$
The {\it slope} of $F$ is defined by
$$\mu_H(F)=\frac{\dg_H(F)}{\rk_H(F)}$$ The sheaf $F$ is {\it
stable} (resp. {\it semistable}) with respect to $H$ if $F$ is
pure of dimension one and for any proper subsheaf
$F'\hookrightarrow F$ one has
$$\mu_H(F')<\mu_H(F) \ (\text{resp.} \leq)$$
For every semistable sheaf $F$ with respect to $H$ there is a {\it
Jordan-H\"{o}lder filtration}
$$0=F_0\subset F_1\subset \hdots \subset F_n=F$$ with stable
quotients $F_i/F_{i-1}$ and $\mu_H(F_i/F_{i-1})=\mu_H(F)$ for
$i=1,\hdots,n$. This filtration need not be unique, but {\it the
graded object} $Gr(F)=\textstyle{\bigoplus_{i}} F_i/F_{i-1}$ does
not depend on the choice of the Jordan-H\"{o}lder filtration. Two
semistable sheaves $F$ and $F'$ on $C$ are said to be {\it
$S$-equivalent} if $Gr(F)\simeq Gr(F')$. Observe that two stable
sheaves are $S$-equivalent only if they are isomorphic. If $F$ is
a semistable sheaf on $C$, we will denote by $[F]$ its
$S$-equivalence class.

By Simpson's work \cite{Si}, there exists a projective moduli
space of semistable pure dimension one sheaves on $C$ of
(polarized) rank 1 and degree 0.

If the curve is not integral, this moduli space can contain some
components given by sheaves of (polarized) rank 1 whose
restrictions to some irreducible components of $C$ are
concentrated sheaves. These components correspond to moduli spaces
of higher rank sheaves on reducible curves (see \cite{L} for some
examples). Therefore, from here on by rank 1 sheaves we mean those
sheaves having rank 1 on every irreducible component of $C$ and we
define {\it the Jacobian of }$C$ as the moduli space
$\overline{M}(C)$ of pure dimension one sheaves of rank 1 and
degree 0 on $C$ that are semistable with respect to the fixed
polarization. Observe that if $C$ is a reduced curve,
$\overline{M}(C)$ is also a projective scheme because it coincides
with Seshadri's compactification \cite{Se}.

For certain projective curves, an explicit description of the
structure of this Jacobian $\overline{M}(C)$ can be found in
\cite{L} where the author studies not only the case of degree 0
sheaves, but also these moduli spaces for arbitrary degree $d$
sheaves.
\subsection{The relative Jacobian} Let $f\colon X\to S$ be an
elliptic fibration and let ${\mathcal O}_X(1)$ be a relative ample
sheaf on $X$. Let $\overline{{\mathcal M}}_{X/S}$ be the functor
which to any $S$-scheme $T$ associates the space of
$S$-equivalence classes of $T$-flat sheaves on $f_T\colon
X\times_S T \to T$ whose restrictions to the fibers of $f_T$ are
semistable of rank 1 and degree 0 with respect to the induced
polarization. Two such sheaves $F$ and $F'$ are said to be
equivalent if $F'\simeq F\otimes f_T^*N$, where $N$ is a line
bundle on $T$.

Again by \cite{Si}, there is a projective scheme
$\overline{M}_{X/S}\to S$ which universally corepresents the
functor $\overline{{\mathcal M}}_{X/S}$. Moreover there is an open
subscheme $\overline{M}^s_{X/S}\subseteq \overline{M}_{X/S}$ that
universally corepresents the subfunctor $\overline{{\mathcal
M}}^s_{X/S}\subseteq \overline{{\mathcal M}}_{X/S}$ of families of
stable sheaves.

The fibration $\hat{f}\colon \overline{M}_{X/S}\to S$ is defined
to be {\it the relative Jacobian} of $f\colon X\to S$. Points of
$\overline{M}_{X/S}$ represent semistable sheaves of rank 1 and
degree 0 on the fibers of $f$ and the natural map $\hat{f}\colon
\overline{M}_{X/S}\to S$ sends a sheaf supported on the fiber
$X_s$ to the corresponding point $s\in S$. In particular, for any
closed point $s\in S$ one has that the fiber $\hat{f}^{-1}(s)$ is
isomorphic to the (absolute) moduli space $\overline{M}(X_s)$ of
semistable rank 1  degree 0 sheaves on $X_s$, that is, the
Jacobian of the curve $X_s$.

In order to study the relative Jacobian $\hat{f}\colon
\overline{M}_{X/S}\to S$ we have then to know the structure of the
corresponding Jacobians of the curves that can appear as fibers of
the elliptic fibration $f\colon X\to S$. We do this in the
following section when $f\colon X\to S$ is an elliptic fibration
of type ($*$) (Definition \ref{d:fibracion}).

\section{The Jacobians of the fibers} Let $f\colon X\to S$ be an
elliptic fibration of type ($*$) and let ${\mathcal O}_X(1)$ be a
line bundle on $X$ ample relative to $S$.

Let $C$ denote a fiber of $f$ and let $H$ be the induced
polarization on $C$. If $C$ is an integral curve (a smooth
elliptic curve, a rational curve with one node or a rational curve
with one cusp), it is well known that the moduli space
$\overline{M}(C)$ is isomorphic to $C$. However if $C$ is a fiber
of type $I_2$, that is, two projective lines meeting transversely
at two points, C\u ald\u araru \cite{Cal2} has proved that the
moduli space $\overline{M}(C)$ is isomorphic to a rational curve
with one node. Following his argument and as a consequence of the
descriptions given in \cite{L}, in this section we prove that for
every reducible fiber $C$ of $f\colon X\to S$, the moduli space
$\overline{M}(C)$ is isomorphic either to a rational curve with
one node or to a rational curve with one cusp.

Let $C$ be any reducible fiber of $f\colon X\to S$, that is, a
curve of type III, IV or I$_N$ with $N\geq 2$. Let us denote by
$C_i$ the irreducible components of $C$. In the following lemma we
collect some properties of rank 1 degree 0 sheaves on $C$ that we
will use later (see \cite{L} for the proof).
\begin{lemma}\label{l:propiedades} If $C$ is a curve of type III,
IV or I$_N$ with $N\geq 2$, it holds that:
\begin{enumerate} \item The (semi)stability of a pure dimension one sheaf of
rank 1 and degree 0 on $C$ does not depend on the polarization.
\item A degree 0 line bundle $L$ on $C$ is stable
if and only if $L|_{C_i}\simeq \mathcal{O}_{\mathbb{P}^1}$ for all
$i$.
\item If $F$ is a stable pure dimension one sheaf of rank 1 and
degree 0 on $C$, then $F$ is a line bundle.
\item If $L$ is a line bundle on $C$ of degree 0, then $L$ is
strictly semistable if and only if $L|_{C_i}\simeq
\mathcal{O}_{\mathbb{P}^1}(r)$ where $r=-1,0$ or $1$ in such a way
that when we remove the components $C_i$ for which $r=0$ there are
neither two consecutive $r=1$ nor two consecutive $r=-1$.
\item If $F$ is a strictly semistable pure dimension one sheaf of
rank 1 and degree 0 on $C$, then its graded object is
$Gr(F)=\oplus_{i}\mathcal{O}_{\mathbb{P}^1}(-1)$
\end{enumerate}
\end{lemma}

Let $q$ be a fixed smooth point of $C$ and let us denote by $C_0$
the irreducible component of $C$ on which $q$ lies.

If $\Delta\subseteq C\times C$ denotes the diagonal and ${\mathcal
J}_\Delta$ is its ideal sheaf, define ${\mathcal O}_{C\times
C}(\Delta)= \Hom ({\mathcal J}_{\Delta},{\mathcal O}_{C\times C})$
as the dual of $\mathcal{J}_\Delta$.

Consider the sheaf $$\mathcal{E}=\mathcal{O}_{C\times
C}(\Delta)\otimes \pi_1^*\mathcal{O}_C(-q)$$ where $\pi_1\colon
C\times C\to C$ is the projection on the first component. This
sheaf is flat over $C$ via the projection $\pi_2\colon C\times
C\to C$ (see \cite{Cal1}, for details) and we have the following

\begin{proposition} \label{p:semiest}
For any point $p\in C$, the restriction $\mathcal{E}_p$ of
$\mathcal{E}$ to $C\times \{p\}$ is a semistable pure dimension
one sheaf of rank 1 and degree 0. Moreover, if $p$ is not a smooth
point of $C_0$, then  all sheaves $\mathcal{E}_p$ define the same
point of the moduli space $\overline{M}(C)$.
\end{proposition}
\begin{proof}
Since $C$ is Gorenstein, we have that
$\ex^1(\mathcal{O}_P,\mathcal{O}_C(-q))=\kappa$, so that the
restriction $\mathcal{E}_p$ is the unique non trivial extension
$$0\to \mathcal{O}_C(-q)\to \mathcal{E}_p\to \mathcal{O}_p\to 0\,
.$$ Using this exact sequence one easily proves that
$\mathcal{E}_p$, which is precisely $\mathcal{J}_p^*\otimes
\mathcal{O}_C(-q)$, is a pure dimension one sheaf of rank 1 and
degree 0.

To prove that it is semistable, let us consider two cases: when
$p$ is a smooth point and when it is a singular point of $C$. In
the first case, $\mathcal{E}_p$ is the line bundle
$\mathcal{O}_C(p-q)$ and we have the following:
\begin{enumerate} \item If $p\in C_0$, the restrictions of $\mathcal{E}_p$ to all irreducible components of $C$ are of
degree 0. Then, by (2) of the previous Lemma, the sheaf
$\mathcal{E}_p$ is stable. \item If $p\notin C_0$, let $C_1$ be
the irreducible component of
 $C$ on which $p$ lies. Since the restriction of $\mathcal{E}_p$
 to $C_0$ has degree -1, to $C_1$ degree 1 and to the others
 components degree 0, (4) in Lemma \ref{l:propiedades} implies that the sheaf
 $\mathcal{E}_p$ is strictly semistable.
\end{enumerate}
In the second case, since $p$ is a singular point of $C$,
$\mathcal{E}_p$ is not an invertible sheaf and then, by (3) in
Lemma \ref{l:propiedades}, it is not stable. Let us see that it is
semistable. Let $\mathcal{G}\hookrightarrow \mathcal{E}_p$ be a
proper subsheaf. We have the exact sequence
$$0\to \mathcal{O}_C(-q)\to \mathcal{E}_p\to \mathcal{O}_p\to 0\,
.$$ Bearing in mind that every proper connected subcurve $D$ of
$C$ has arithmetic genus 0 and $D\cdot \overline{D}=2$, it is easy
to deduce from Lemma 3.4 in \cite{L} that the line bundle
$\mathcal{O}_C(-q)$ is stable. Consider the composition map
$g\colon \mathcal{G}\to \mathcal{O}_p$ that can be either zero or
surjective. If $g$ is zero, $\mathcal{G}\hookrightarrow
\mathcal{O}_C(-q)$ and then $\mu_H(\mathcal{G})<-1$. If $g$ is
surjective and we denote by $\mathcal{H}$ its kernel, we have that
$\mathcal{H}$ is a subsheaf of $\mathcal{O}_C(-q)$ and then
$\mu_H(\mathcal{H})<-1$. Since the degree of $\mathcal{H}$ is
integer, from the exact sequence
$$0\to \mathcal{H}\to \mathcal{G}\to \mathcal{O}_p\to 0\, ,$$ we
conclude that $\mu_H(\mathcal{G})\leq 0$. Then  $\mathcal{E}_p$ is
a strictly semistable sheaf.

For the second part of the statement, it is enough to note that if
$p$ is not a smooth point of $C_0$, $\mathcal{E}_p$ is a strictly
semistable sheaf and then, by (5) in \ref{l:propiedades}, its
graded object is isomorphic to
$\oplus_i\mathcal{O}_{\mathbb{P}^1}(-1)$. Hence all these sheaves
are in the same $S$-equivalence class and the proof is complete.
\end{proof}

The restriction of the family $\mathcal{E}$ to $C\times C_0$
gives, by the universal property of $\overline{M}(C)$, a map
$$\phi\colon C_0\to \overline{M}(C)$$ defined as
$\phi(p)=[\mathcal{E}_p]$.

The same proof that C\u ald\u araru gives in \cite{Cal2} when $C$
is a curve of type $I_2$ proves the following
\begin{proposition}If $C$ is any reducible fiber of a elliptic
fibration of type ($*$), then the sheaves $\mathcal{E}_p$ satisfy
\begin{equation}
\ex^i(\mathcal{E}_p,\mathcal{E}_{p'})=
\begin{cases}
\kappa \text{ \ \ if } p=p' \text{ and } i=0  \\
0 \text{ \ \ if } p \neq p' \text{ and all } i\notag
\end{cases}
\end{equation} In particular the
moduli space $\overline{M}(C)$ has dimension 1 and is smooth at
$[\mathcal{E}_p]$ for any smooth point $p\in C_0$.
\end{proposition}

Thus if $p$ and $p'$ are two different smooth points of $C_0$,
since the sheaves $\mathcal{E}_p$ and $\mathcal{E}_{p'}$ are
stable and $\ho(\mathcal{E}_p,\mathcal{E}_{p'})=0$, we have that
$\phi(p)\neq \phi(p')$. Then the map $\phi$ is injective in
$C_0\setminus \overline{C_0}$ where $\overline{C_0}$ denotes the
complementary subcurve of $C_0$ in $C$. Since $C_0$ is irreducible
and $[\mathcal{O}_C]\in \im \phi$, the map $\phi$ factors as
$$\xymatrix{ C_0\ar[r]\ar[d] & \overline{M}(C)\\
M'(C)\ar@{^{(}->}[ur]}$$ where $M'(C)$ is an irreducible component
of $\overline{M}(C)$ that contains the point $[\mathcal{O}_C]$.
But since $M'(C)$ is of dimension 1 and smooth at
$[\mathcal{O}_C]$, it is the unique irreducible component of
$\overline{M}(C)$ containing $[\mathcal{O}_C]$ and the map
$\phi\colon C_0\to M'(C)$ is also surjective.

In this point, we distinguish to cases:

1- If $C$ is a curve of type I$_N$ with $N\geq 2$, let $\{r_1,
r_2\}$ be the two intersection points between $C_0$ and
$\overline{C_0}$. By Proposition \ref{p:semiest}, we have that
$\phi(r_1)=\phi(r_2)$ and then $M'(C)$ is a rational curve with
one node (Figure 1).

\hspace{1truecm}
\centerline{\includegraphics[scale=0.3]{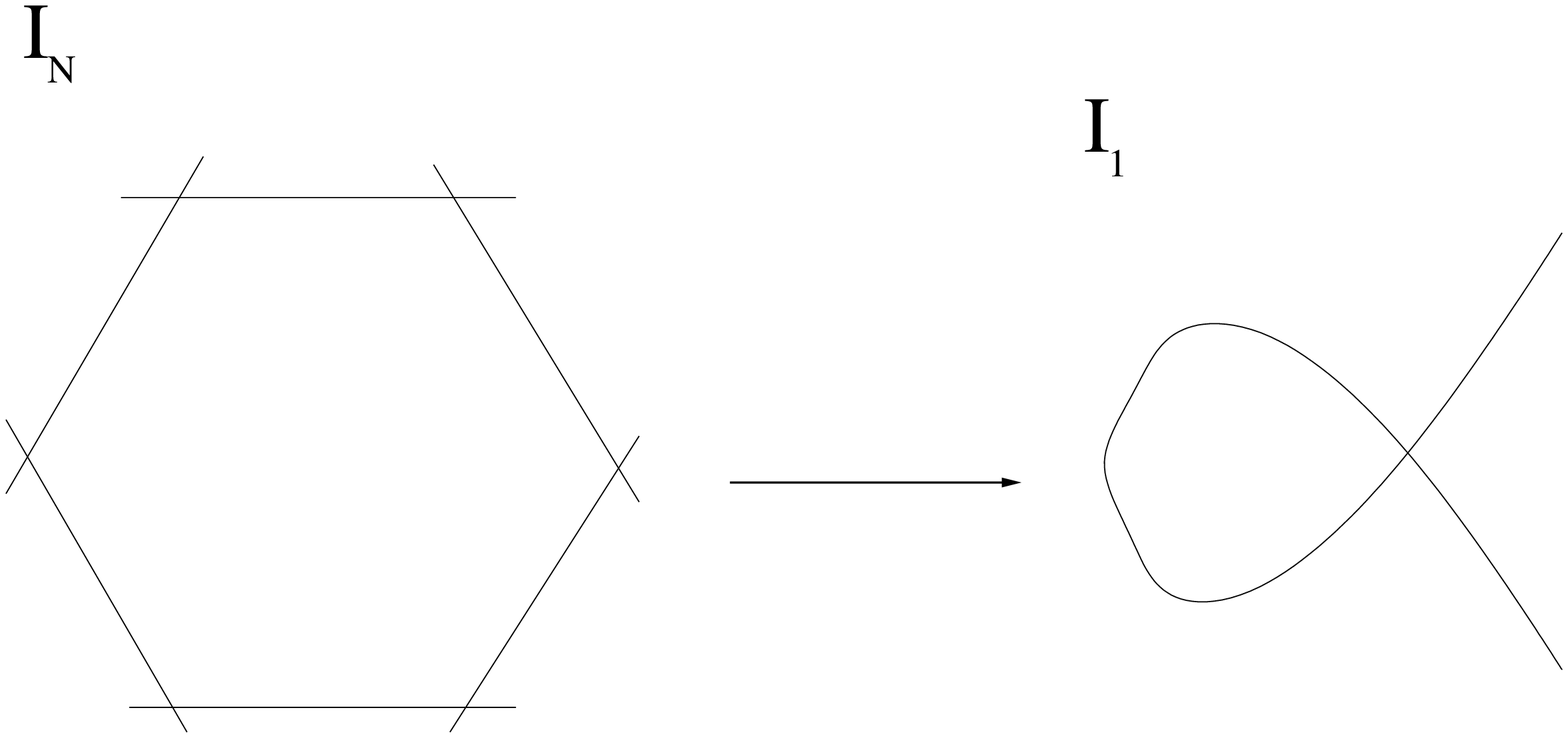}\hspace{3truecm}}

\hspace{0.5truecm}

 \centerline{   Figure 1.}

\hspace{1truecm}

 In fact, since by Proposition 5.13 in \cite{L}, stable line
bundles on $C$ of degree 0 are given by the group exact sequence
$$0\to \mathbb{G}_m\to \overline{M}^s(C)\to \prod_{i=1}^N
\Pic^0(C_i)\to 0$$ and, by Corollary 6.7 in \cite{L}, there is
only one extra point in $\overline{M}(C)$ corresponding to any
strictly semistable sheaf, we can conclude that
$\overline{M}(C)=M'(C)$. Thus the moduli space $\overline{M}(C)$
is in this case isomorphic to a rational curve with one node.

2- If $C$ is a curve of type III or IV and $r$ is the intersection
point of $C_0$ and $\overline{C_0}$, we know that $\phi(r)$ is the
unique singular point of $M'(C)$ and that it is a cusp (Figure 2).

\hspace{2truecm}

\centerline{\includegraphics[scale=0.3]{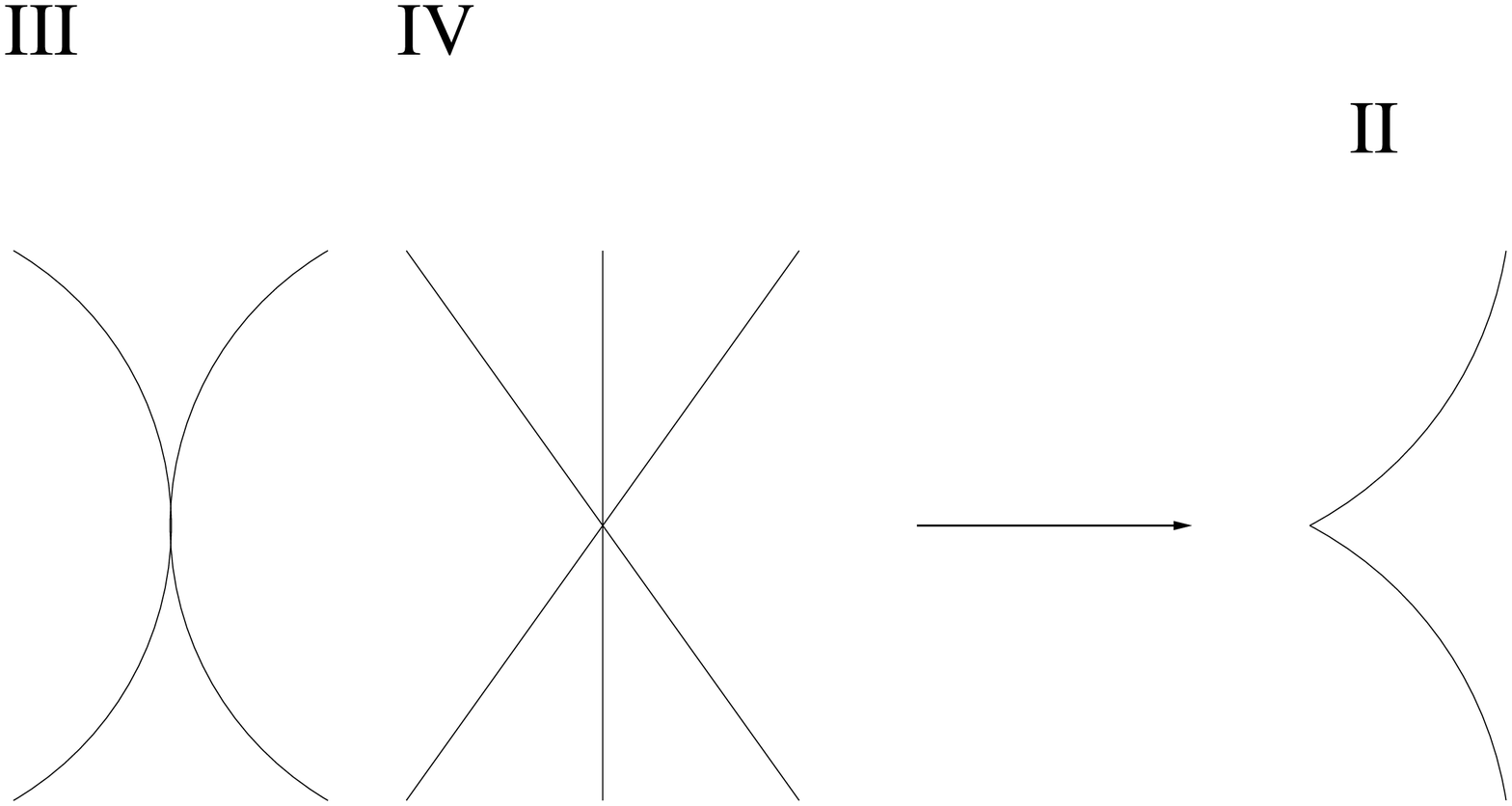}\hspace{0.5truecm}}

\centerline{   Figure 2.}

\hspace{1truecm}

Moreover, Propositions 5.6 and 5.8 of \cite{L} imply that stable
line bundles on $C$ of degree 0 are determined by the group exact
sequence:
$$0\to \mathbb{G}_a\to \overline{M}^s(C)\to \prod_{i=1}^N
\Pic^0(C_i)\to 0$$ and as above there is only one extra point in
$\overline{M}(C)$. Then we conclude that $\overline{M}(C)=M'(C)$
is isomorphic to a cuspidal irreducible curve.

In this case we don't know a priori the arithmetic genus of
$\overline{M}(C)$. However in the next section we will prove that
the relative Jacobian $\hat{f}\colon \overline{M}_{X/S}\to S$ of
$f$ is a flat morphism, and then all its fibers have arithmetic
genus 1. This allows us to conclude  that the moduli space
$\overline{M}(C)$ is isomorphic to a rational curve with one cusp.

\section{The Flatness of the relative Jacobian} Let $f\colon X\to
S$ be an elliptic fibration of type ($*$) and $\hat{f}\colon
\overline{M}_{X/S}\to S$ its relative Jacobian. From the previous
section we know that this relative Jacobian is a projective
morphism whose fibers are integral curves. We are now going to
prove that it is a flat morphism.

Using Corollary 15.2.3 in \cite{EGA} about the flatness of
universally open morphisms with reduced fibers, in order to prove
that $\hat{f}$ is flat, we only need to see that it is universally
open. Taking into account that $S$ is a smooth variety and then
geometrically unibranch, if $\hat{f}$ is equidimensional we
conclude thanks to Chevalley's criterion (Corollary 14.4.4 in
\cite{EGA}). Since $\hat{f}$ is surjective and all its fibers have
the same dimension, to show that it is equidimensional it is
enough to prove that $\overline{M}_{X/S}$ is irreducible.

Notice first that, since for all $s\in S$ the sheaf
$\mathcal{O}_{X_s}$ is stable, the flat family $\mathcal{O}_X$
defines a natural section $\sigma\colon S\hookrightarrow
\overline{M}_{X/S}$  of $\hat{f}\colon \overline{M}_{X/S}\to S$
whose image is contained in only one irreducible component $M'$ of
$\overline{M}_{X/S}$. Indeed, since $S$ is irreducible,
$\im\sigma$ is contained in some irreducible component of $
\overline{M}_{X/S}$. But for all $s\in S$,
$\sigma(s)=[\mathcal{O}_{X_s}]$ is a smooth point of
$\overline{M}_{X/S}$ so that this irreducible component $M'$ is
unique.

Let us see that $\overline{M}_{X/S}=M'$. Actually, since
$(\overline{M}_{X/S})_s$ is irreducible for every $s\in S$, it is
contained in some irreducible component of $\overline{M}_{X/S}$,
namely $\tilde{M}$. If $\tilde{M}\neq M'$, the point
$[\mathcal{O}_{X_s}]$, which lies on $(\overline{M}_{X/S})_s$ and
on the image of the section $\sigma$, belongs to $\tilde{M}\cap
M'$. But this is not possible because this is a smooth point of
$\overline{M}_{X/S}$. Hence all fibers $(\overline{M}_{X/S})_s$
are contained in the irreducible component $M'$ and so it is the
whole moduli space $\overline{M}_{X/S}$.

The final result is then the following
\begin{theorem}\label{teorema} If $f\colon X\to S$ is an elliptic fibration of
type ($*$) (Definition \ref{d:fibracion}), the moduli space of
semistable pure dimension one sheaves of rank 1 and degree 0 on a
fiber $X_s$ is isomorphic to:
\begin{enumerate}
\item A smooth elliptic curve if $X_s$ is smooth.
\item A rational curve with one node if $X_s$ is of type I$_N$
with $N\geq 1$.
\item A rational curve with one cusp if $X_s$ is of type II, III,
or IV.
\end{enumerate}
\end{theorem}

Thus the relative Jacobian $\hat{f}\colon \overline{M}_{X/S}\to S$
of an elliptic fibration of type ($*$) is an integral elliptic
fibration which always has a global section even if the original
fibration has no sections. Here the difference between the
integral case and the case with reducible fibers is the following.
If the integral elliptic fibration has a global section, we know
that it is globally isomorphic to its relative Jacobian, in
contrast if the original fibration has reducible fibers even if it
has a section, we are not able to ensure that it is globally
isomorphic to its relative Jacobian because, as we will see now,
$\overline{M}_{X/S}$ can be a singular space.

In fact, this theorem shows that as long as the original fibration
$f\colon X\to S$ has reducible fibres, to get its relative
Jacobian $\hat{f}\colon \overline{M}_{X/S}\to S$  we have to
contract to a point all but one irreducible components of every
reducible fibre $X_s$ of $f$, that is, for every reducible fibre
of $f$ we have to contract to a point $q\in\overline{M}_{X/S}$ a
linear chain $\cup_i C_i$ of smooth rational curves.

When $f\colon X\to S$ is a smooth elliptic surface, being
$C_i^2=-2$ for every $i$, we know  by \cite{A} that the
contraction of a such chain is a singular point. Since the
discriminant locus of $f$ is a finite number of points, we
conclude that the relative Jacobian $\hat{f}\colon
\overline{M}_{X/S}\to S$ is an integral elliptic surface with at
worst a finite number of singular points. When $f\colon X\to S$ is
a smooth elliptic threefold, by \cite{M} every irreducible
component of these chains has length 1, and then using the results
in \cite{Z} we have that the contraction point $q$ is also a
singular point. However in this case, since the discriminant locus
of $f$ has dimension one, the singular locus of the relative
Jacobian $\overline{M}_{X/S}$ has dimension less or equal to 1.

\end{document}